\documentclass[letterpaper, 10 pt, conference]{ieeeconf}  

\IEEEoverridecommandlockouts                              

\usepackage{amsfonts}
\usepackage{amssymb,amsmath}
\usepackage{graphicx}
\usepackage{cite}
\renewcommand{\epsilon}{\varepsilon}

\renewcommand{\d}{\,\mathrm{d}}
\newcommand{\id}{\mathop{\bf id}}

\newtheorem{theorem}{Theorem}
\newtheorem{remark}{Remark}

\renewcommand{\epsilon}{\varepsilon}

\newcommand{\spt}{\mathop{\rm spt}}

\renewcommand{\d}{\,d}

\newcommand{\R}{\mathbb{R}}

\newcommand{\C}{\mathcal{C}}

\title{\LARGE \bf
Optimal control of distributed ensembles\\ with application to Bloch equations$^*$
}

\author{Roman Chertovskih$^{1}$, Nikolay Pogodaev$^{2}$, Maxim Staritsyn$^{1}$, and A. Pedro Aguiar$^{1}$
\thanks{*The authors acknowledge the financial support of the Foundation for Science and
Technology (FCT, Portugal) 
in the framework of the Associated Laboratory ``Advanced Production and Intelligent
Systems'' (AL ARISE, ref. LA/P/0112/2020), R\&D Unit SYSTEC (base UIDB/00147/2020 and programmatic UIDP/00147/2020
funds), and projects SNAP (ref. NORTE-01-0145-FEDER-000085) and MLDLCOV (ref. DSAIPA/CS/0086/2020).
}
\thanks{$^{1}$Roman Chertovskih, Maxim Staritsyn, and A. Pedro Aguiar are with Research Center for Systems and Technologies (SYSTEC), ARISE \& Department of Electrical
and Computer Engineering,
Faculdade de Engenharia, Universidade do Porto, Rua Dr. Roberto Frias, s/n 4200-465, Porto, Portugal
        {\tt\small roman@fe.up.pt, staritsyn@fe.up.pt, pedro.aguiar@fe.up.pt}
}%
\thanks{$^{2}$Nikolay Pogodaev is with Dipartimento di Matematica ``Tullio Levi-Civita'' (DM), University of Padova, Via Trieste, 63 - 35121 Padova, Italy
        {\tt\small nickpogo@gmail.com}}%
}

\begin{document}
\maketitle
\thispagestyle{empty}
\pagestyle{empty}
\begin{abstract}
Motivated by the problem of designing robust composite pulses for Bloch equations in the presence of natural  perturbations, we study an abstract optimal ensemble control problem in a probabilistic setting with a general nonlinear performance criterion. The model under study addresses mean-field dynamics described by a linear continuity equation in the space of probability measures. For the resulting optimization problem, we derive an exact representation of the increment of the cost functional in terms of the flow of the driving vector field. Relying on the exact increment formula, a descent method is designed that is free of any internal line search. 
The numerical method is applied to solve new control problems for distributed ensembles of Bloch equations. 
\end{abstract}
\section{MOTIVATION}

Consider a 
population of homotypic individuals labeled by the points $\omega$ of some set $\Omega$. The state of the $\omega$th object at the time moment $t$, $x(t,\omega) \in \mathbb R^n$, evaluates on a given time interval $I \doteq [0, T]$ under the action of the parameterized vector field
$V: \, \R^n \times 
\R^s \times U \to \R^n$, starting from a given position $x_0(\omega) \in \R^n$: 
\begin{equation}
\left\{\begin{array}{l}
\partial_t x(t,\omega) = V_{u}\left(x(t,\omega),\eta(\omega)\right) \\[0.2cm] x(0,\omega) = x_0(\omega)
\end{array}\right| \quad \omega \in \Omega.\label{rand-ODE}
\end{equation}
The dynamics \eqref{rand-ODE} involves two types of structural ``parameters'': 
the function $\eta: \, \Omega \to \R^s$ manifests  disturbances or structural variations of the underlying model, while 
an \emph{exogenous} signal $u$ with values in a given set $U \subseteq \R^m$ models the \emph{control} action. 

In the simplest case, the parameterization space $\Omega$ is just a finite set of indexes, and \eqref{rand-ODE} reduces to a multi-agent system of non-interacting units. In a more general setup, we deal with 
the \emph{continuum} of 
individuals moving in a \emph{dis}coordinated way. 
Commonly, in such models, $\Omega$ is a simply organized compact subset of $\R^s$,
and $\eta$ is the identity mapping $\Omega \to \Omega$. 

Problems of \emph{ensemble control} arise when one has to design a control signal in a ``broadcast'' way, i.e.,
such that it acts simultaneously on all individual trajectories $x(\cdot,\omega)$, $\omega \in \Omega$, 
to force them towards a desired behavior; 
this means that $u$ should be a function $t \mapsto u(t)$ of time variable only (independent of $\omega$).

A canonical example is the problem of designing external excitations of quantum ensembles. Pioneering works in this area were focused on the famous Bloch equation~\cite{Li-Khaneja,Li-Khaneja2006}, which models the macroscopic evolution of bulk magnetization in a population of non-interacting nuclear spins immersed in an intense static magnetic field, which is modulated by the radio frequency (rf-) field. In 
nuclear magnetic resonance (NMR) experiments, the strength of the applied magnetic field is subject to unavoidable perturbations (static- and/or rf-field inhomogeneity), while the spin ensembles demonstrate perceptible variations in their dissipation rates and/or natural frequencies (Larmor dispersion). The related problem of control engineering is to design \emph{robust} signals (so-called composite pulses) compensating for the mentioned disturbances; mathematically, this task can be formalized as a problem of optimal ensemble control, 
see, e.g. \cite{SKINNER20038}. 
In NMR spectroscopy, the designed pulse sequences are typically desired to be selective, i.e., some sub-populations (with prescribed Larmor frequencies) have to be excited, while the other ones should remain intact or saturated \cite{Conolly}; such are, e.g., contrast problems in NMR imaging  \cite{Bonnard-TAC,Bonnard2016}. In the language of ensemble control, this means to drive several \emph{uncoupled} populations of spins by a common magnetic field.

\subsection{Probabilistic Setup. 
Distributed Ensembles}\label{sec:intro-1} 

In contrast to \cite{WANG2018306,Bonnard2016,Bonnard-TAC}, 
our approach stems from the {probabilistic interpretation} of the ensemble dynamics, assuming that $\Omega$ is endowed with the structure of probability space $(\Omega, \mathcal A, \mathbb P)$ with a specified $\sigma$-algebra $\mathcal A \subset 2^\Omega$ and a canonical probability measure $\mathbb P$ on $(\Omega, \mathcal F)$ (we shall write $\mathbb P \in \mathcal P(\Omega)$). 

This interpretation is motivated by practical applications, in which the individual states $x(\cdot, \omega)$ can not be measured directly, 
and
all the available information
is based on some ``observables''~-- measurement outputs accompanying the dynamics \eqref{rand-ODE}
and involving 
certain statistical characteristics, see, e.g., 
\cite{CHEN2020109057}. 

In the probabilistic setup, the map \((t, \omega) \mapsto (x(t,\omega), \eta(\omega))\) is naturally 
viewed 
as a deterministic random process, and the behavior of the random variable $\omega \mapsto (x(t,\omega),\eta(\omega))$ can be analyzed by investigating the time-evolution of its law 
\begin{align}\label{mu-Law}
\varrho_t = \left(x(t,\cdot), \eta(\cdot)\right)_{\sharp} \mathbb P \in \mathcal P(\mathbb R^{n+s}).
\end{align} 
Hereinafter, the operator $F_\sharp: \, \mathcal P(\mathcal X) \to \mathcal P(\mathcal Y)$ 
denotes the \emph{pushforward} of a 
measure $\mu\in \mathcal P(\mathcal X)$ 
through a (Borel) map $F :\, \mathcal X \to \mathcal Y$ between two measurable spaces that acts on functions $\varphi: \, \mathcal Y \to \R$ with the property $\varphi\circ F \in L_1(\mathcal X;\mu)$ by the rule 
\begin{equation}
\int_{\mathcal Y} \!\! \varphi \d (F_\sharp \mu) = \int_{\mathcal X}\!\! \varphi \circ F \d \mu.\label{sharp}
\end{equation}
Under the standard regularity of the map $(x,\eta) \mapsto V_\upsilon(x,\eta)$, the measure-valued curve $t \mapsto \varrho_t$ is a unique distributional solution of the 
\emph{continuity equation} \cite{AmbrosioBook2007}
\begin{equation}
    \partial_t \varrho_t + \nabla_x \cdot \left( V_{u(t)} \, \varrho_t \right) =0, \quad \varrho_0  = (x_0(\cdot),\eta(\cdot))_{\sharp}\mathbb P;\label{PDE}
\end{equation}
$\nabla_x$ denotes the gradient w.r.t. $x \in \mathbb R^n$, and ``$\cdot$'' means the scalar product. 

The discussed interpretation of \eqref{rand-ODE} postulates a passage from the multi-particle, \emph{microscopic} model represented by many copies of an ODE to a  distributed, \emph{macroscopic} representation described by  a PDE and called the \emph{mean field}; systems \eqref{rand-ODE} and \eqref{PDE} are the so-called Lagrangian and Eulerian forms of the mean-field dynamics, respectively \cite{Cavagnari2018}. 

Remark that, as a result of this passage, a nonlinear finite-dimensional object is replaced  by an infinite-dimensional but  \emph{state-linear} ($\varrho$-linear) one. The linearity of the reduced model plays a vital role in our study as it gives rise to an exact representation of the increment (\emph{$\infty$-order} variation) of the cost functional in the corresponding optimal control problem to be presented in \S~\ref{sec:formula}.  

Finally, observe that PDE \eqref{PDE} can be viewed as a family $\eta \mapsto \mu^\eta_t$ of $ \mathcal P(\R^n)$-valued curves solving, $\Xi$-a.e., the ``sliced'' continuity equation of the same structure with the vector field $V^\eta \doteq V(\cdot, \eta)$ and  initial condition $\mu^\eta_0 = \vartheta^\eta \in \mathcal{P}(\R^n)$, where the map $\eta \mapsto \vartheta^\eta$ is obtained by disintegrating the distribution $\varrho_0$ w.r.t. the projection $\Xi \doteq ((x, \eta) \mapsto \eta)_\sharp \varrho_0$. We call such a family the \emph{distributed ensemble}; this concept separates two types of uncertainty: dispersion in the initial data $\omega \mapsto x_0(\omega)$ is converted to the mean field, while fluctuations of the dynamics, $\eta \mapsto V^\eta$, are treated independently. 

\subsection{Contribution and Novelty}

This work contributes to the line of research \cite{WANG2018306,Bonnard-TAC,Bonnard2016,SKINNER20038} devoted to optimal control of quantum ensembles. We elaborate on a general approach  that captures the natural probabilistic flavor of ensemble control problems. 
From practical viewpoints, it enables us to improve the quality of designed control signals since it takes into account the available statistical information, and in this way allows us 
to concentrate the ``resource'' of feasible control options around \emph{relevant} values of $\omega$. A key result is the development of a descent algorithm for optimal ensemble control originating from an exact increment formula for the \emph{nonlinear} cost functional. In contrast to familiar indirect methods \cite{Bonnet2021AMT} based on the 1st variation (i.e. on Pontryagin's maximum principle, PMP), our approach is free of any hidden parameters and does not involve any internal line search. This essentially improves the computational performance, where the algorithm is proved to converge towards a PMP extremal, but the convergence is, typically, faster than as for the conventional gradient descent. Furthermore -- due to the nonlocal nature of the underlying increment formula -- our algorithm can step over local solutions, and therefore, has the potential of global search.\footnote{Since the formula is exact, the generated control variations should \emph{not} be sufficiently ``small''; they also should \emph{not} be of any specific class such as needle-shaped or weak variations, as it is common for the classical optimal control theory.} 

This paper generalizes our recent works \cite{SChPP-2022,SPP-2023}, where the exact increment formula and a nonlocal algorithm were derived for models of linear and linear-quadratic structure. Now, we consider 
an arbitrary nonlinear cost functional on the space of probability measures, which 
has the so-called intrinsic derivative (see \cite{CardMaster2019} and the discussion in sec.~\ref{sec:diff-meas}). 

\section{OPTIMAL CONTROL PROBLEM}

First, we introduce some necessary \underline{notations}: Let  $\mathcal X$ be a 
metric space, and $I \doteq [0,T]$. We denote by $\C(I; \mathcal X)$ the spaces of continuous maps $I\mapsto \mathcal X$ with the usual $\sup$-norm. If $\mathcal X \subseteq \R^n$, $\C^1(\mathcal X)$ denotes the space of continuously differentiable functions $\mathcal X \to \R$, and  $\C^\infty_c(\mathcal X)$ the space of smooth functions with a compact support in $\mathcal X$; $L_p(I;\mathbb R^m)$, $p=1,\infty$, the Lebesgue 
spaces of summable and bounded measurable functions $I\mapsto \mathbb R^m$, respectively. 

$\mathcal P(\mathcal X)$ the set of probability measures on $\mathcal X$, and $\mathcal P_c(\mathcal X) \subseteq \mathcal P(\mathcal X)$ the set of measures having compact support in $\mathcal X$; $\mathcal P_c(\R^n)$ is a complete separable metric space as it is endowed with any $p$-Kantorovich (Wasserstein) distance $W_p$, $p \geq 1$. 

Among all measures on $\R^n$, we mark out two specific ones~-- the usual Lebesgue measure, $\mathcal L^n$, and a Dirac point-mass measure  concentrated at $x \in \R^n$, $\delta_x$. 

\subsection{General Problem Statement}

Our prototypic mathematical object is the following optimization problem on $\mathcal P_c(\R^n)$: 
\begin{align}
(P) \quad \min \quad & \mathcal I [u]=\ell(\mu_T) \nonumber\\
\mbox{subject to }\quad &\partial_t\mu_t + \nabla_x\cdot \left(V_{u}\, \mu_t\right) = 0, \label{conteq}\\
&t \in I\doteq [0,T]; \quad \mu_0=\vartheta;\label{mu-init}\\
& u(\cdot) \in \mathcal U \doteq L_\infty(I; U), \ U\subset \R^m,
\end{align}
where $\ell: \, \mathcal P_c(\R^n) \to \R$ is a given performance criterion, and $V: \, \R^n \times U \to \R^n$ a control vector field.\footnote{With slight abuse of notation, we use the letter $V$ in different contexts.} 
Despite its probabilistic appearance, $(P)$ is a \emph{deterministic} optimal control problem, in which 
the trajectories are measure-valued functions $\mu \in C(I, \mathcal P_c(\R^n))$, and the control signals are usual functions $u \in L_\infty(I, U)$. 
This problem 
can be specified to the case of distributed ensembles as follows: 
\begin{equation}
(\widetilde P) \quad \min\int_{\R^s} \ell(\mu^\eta_T) \d \Xi(\eta),\label{cost-exp}
\end{equation}
where $t \mapsto \mu_t^\eta[u]$ solves the linear PDE \eqref{conteq}, \eqref{mu-init} with $V_u=V^\eta_u$ for $\Xi$-a.a. $\eta \in \R^s$.  

We make the following standard \underline{regularity hypotheses}:
\begin{itemize}
    \item[$(A_1)$] the map 
    $(x, \upsilon)\mapsto V_\upsilon(x)$ is continuous, continuously differentiable in $x$ and satisfies the sublinear growth condition: there exists a constant $M>0$ such that
    \(
    V_\upsilon(x) \leq M(1+|x|)
    \) 
    for all $x \in \R^n$ and $\upsilon \in U$.  
\item[$(A_2)$] The set $U$ is convex and compact.

\item[$(A_3)$] $\vartheta  \in \mathcal P_c(\R^n)$, and $\ell: \, \mathcal P_c(\R^n) \to \R$ is continuous.

\item[$(A_4)$] $\ell \in \C^1$ in the sense of intrinsic derivative (to be specified below).
\end{itemize}

$(A_1)$ is the standard 
set of assumptions to guarantee the well-posedness of the PDE \eqref{conteq} \cite{AmbrosioBook2007}. $(A_1)$--$(A_3)$ imply the existence of a minimizer for problem $(P)$  \cite[Theorem~3.2]{POGODAEV20203585}; 
under these assumptions, the solution $\mu_t$ of \eqref{conteq}, \eqref{mu-init} is supported in a ball whose radius depends only on the problem data  \cite[Lemma~A2]{POGODAEV20203585}. Hence, $\mu_t \in \mathcal P_c(\R^n)$ for all $t \in I$.

\subsection{Problem Specification}\label{sec:diff-problems}

Below, we provide some examples of the performance criterion $\ell$ that cover typical optimization tasks in the area of ensemble control. 

\subsubsection*{Targeting}

In several NMR applications, the guide is supposed to transfer the ensemble from one given profile $x(0,\cdot) = x_0(\cdot)$  (as close as possible) to another one $x(T,\cdot) = x_T(\cdot)$, which means to minimize the quantity
\[
\displaystyle\int_\Omega |x(T,\omega) - x_T(\omega)|^2 \d \mathbb P(\omega).
\]
Such are problems of selective spin excitation, see \cite{Li-TACON-2013,WANG2018306} and the bibliography therein.

This problem is formulated in our setting by using \eqref{mu-Law}, the definition of $\mu^\eta$, and the change of variable formula \eqref{sharp}:
\begin{equation}
\min\int_{\R^s} \ell(\mu^\eta_T; x_T(\eta))\d \Xi(\eta)\label{type-1},
\end{equation}
where
\(
\ell(\mu;x) \doteq  \displaystyle\int_{\R^n}|y-x|^2\d \mu(y)\).
Typically, the map $x_T$ is chosen to be constant, which means that the ensemble is assumed to be aggregated around some given position.

\subsubsection*{Statistical Tracking}

In some cases \cite{SPP-2023,ZUAZUA20143077}, the previous performance criterion could be too rigid. 
Instead of matching  the desired profile in average, one may require that the target distribution has prescribed statistical characteristics, for instance, its expectation and variance approach some desired values. The cost functional can be reset in the language of distributed ensembles as follows:  
\begin{equation}
\int_{\R^s} \!\! \left[\psi_1\left(\mathcal E(\mu_T^\eta)- \hat{\mathcal E}\right) + \psi_2\left(\mathcal V(\mu_T^\eta)- \hat{\mathcal V}\right)\right]\d \Xi(\eta),\label{type2-problem}
\end{equation}
where $\mathcal E(\mu)$ and $\mathcal V(\mu)$ denote the expectation and variance of $\mu \in \mathcal P(\R^n)$, respectively,
$\hat{\mathcal E} \in \R^n$ and $\hat{\mathcal {V}} \in \R$ are target values of the statistical characteristics, and
$\psi_{1}: \, \R^n \to \R$ and $\psi_2: \, \R \to\R$ are given penalty functions. 

\subsubsection*{Minimum-Energy Control} In many applications, the discussed cost functionals are accompanied by the energy term
\begin{align}
    \frac{\alpha}{2}\int_0^T |u(t)|^2 \d t\label{u-cost}
\end{align}
with some weight $\alpha>0$. In particular, this produces a sort of \emph{regularization} of the underlying 
problem.

\section{PRELIMINARIES}

In this section, we provide the necessary theoretical background and collect some auxiliary results.

\subsection{Flows of Vector Fields. Transport Equation}

Let $V: \, I \times \R^n \to \R^n$ be a time-dependent vector field generating a flow, i.e. a map $X: \, I\times I \times \R^n \to \R^n$ such that, for all $s \in \R$ and $x \in \R^n$, the function $t \mapsto X_{s,t}(x)$ is a solution of the ODE
\begin{equation}
\partial_t X_{s,t} = V_t\circ X_{s,t}, \quad X_{s,s}=\id,\label{flow}
\end{equation}
where $\id$ stands for the identical map $\R^n \to \R^n$. In view of the semigroup property 
\(
X_{t_0,t_2} = X_{t_1,t_2}\circ X_{t_0,t_1}\)  \(\forall \, t_0, t_1, t_2,
\)
the inverse of $X_{s,t}$ is the map $X_{t,s}$. 

Fixed $s$, abbreviate $P_t =X_{s,t}$ and $Q_t = X_{t,s}$. Then, by the chain rule,
\(
0=\partial_t(\id) = \partial_t(Q_t \circ P_{t}) = (\partial_t Q_t + D_x Q_t\, V_t)\circ P_{t}.\) 
Since the expression in the brackets vanishes  for all values $P_{t}(x)$, and therefore, for any $x \in \R^n$, we conclude that the inverse flow should satisfy the linear operator equation
\begin{equation}
\partial_t Q_t + D_x Q_t\, V_t =0, \quad Q_s=\id.\label{dual-transp-loc}
\end{equation}
Returning to the $X$-notation,
and recalling that the Jacobian $J_{t,s} \doteq D_x X_{t,s}$ satisfies \cite[Ths. 2.2.3 and 2.3.2]{Bressan-PiccoliBOOK} the linear problem
\begin{equation}
\partial_t J_{t,s} = -  J_{t,s}\, (D_x V_t \circ X_{t,s}), \quad J_{s,s}=E,\label{Jacobian}
\end{equation}
where $E$ denotes the identity matrix, we express the derivative of the inverse flow w.r.t. $t$ as follows:
\begin{equation}\label{inv-flow-derivative}
\partial_t X_{t,s} = - J_{t,s}\, V_t.
\end{equation}

Note that operators $P$ and $Q$ refer to the concepts of the left and right chronological exponents in the tradition of geometric control theory \cite{Agrachev-Book}. 

\subsection{Continuity Equation}\label{Sec-dual}

Recall that the continuity equation \eqref{conteq} on the space $\mathcal P_c(\R^n)$ is understood in the weak (distributional) sense. A function $\mu: \, t \to \mu_t$ is said to be a weak solution of  \eqref{conteq} iff the following equality holds 
\begin{equation}
\label{conteq-sol}
\begin{array}{c}\displaystyle 0 = \int_0^T  \d t \int \left(\partial_t \varphi_t + \nabla_x \varphi_t \cdot V_{u(t)}\right)\d \mu_t 
\end{array}
\end{equation}
for all $\varphi \in C^\infty_c((0,T)\times \R^n)$; hereinafter, we abbreviate $\int = \int_{\R^n}$. 
Under 
assumptions $(A_1)$, there exists a unique weak solution to \eqref{conteq} with initial condition \eqref{mu-init};
this solution admits the following representation  \cite{AmbrosioBook2007} in terms of the characteristic flow \eqref{flow}:
\(
\mu_t = (X_t)_\sharp \vartheta,
\) where 
$X_t \doteq X_{0,t}$.

\subsection{Differentiation w.r.t. the Probability Measure}\label{sec:diff-meas}

Since $\mathcal P_c$ is merely a metric space and does not have a linear structure, standard concepts of the directional derivative are not applicable here (there are simply no ``directions'' in common sense). 
At the same time, there is an option to differentiate a function $F: \, \mathcal P_c(\R^n) \to \R$ at some $\mu \in \mathcal P_c(\R^n)$ in the ``direction'' of a (Borel measurable and locally bounded) vector field  \( f\colon \R^{n}\to \R^{n} \)  pushing the measure $\mu$: 
$\frac{d}{d\lambda}\Big|_{\lambda=0}F\left((\id+\lambda f)_{\sharp}\mu\right)$.
Under some reasonable regularity \cite{CardMaster2019} of the map $F$, this  derivative does exist and takes the form: 
$
\displaystyle\int D_{\mu}F(\mu)\cdot f\d\mu,$
where the linear map \( D_{\mu}F\colon \mathcal{P}_{c}(\R^{n})\times\R^{n}\to\R^{n} \), called the \emph{intrinsic derivative}, can be calculated as follows
\begin{align}
D_{\mu} F(\mu)(x)  = 
D_{x}\lim_{h\downarrow 0}\frac{1}{h}\big(F\left(\mu\!+\!h(\delta_{x}\!-\!\mu)\right)\!-\!F(\mu)\big).\label{intrins-der}
\end{align}
The expression under the sign of $D_x$ is called the \emph{flat derivative} of $F$ (typically denoted by \(\frac{\delta F}{\delta\mu}\)). Note that the notions of intrinsic and flat derivatives are naturally connected to another useful concept of derivative on $\mathcal P_c(\R^n)$, the so-called \emph{localized Wasserstein derivative}~\cite{Bonnet2021}.  

In contrast to the other concepts of derivative in the space of measures, the quantity \eqref{intrins-der} can be computed explicitly (and rather easily) for many functionals arising in practice, in particular, for those specified in \S~\ref{sec:diff-problems}. Below, we shall utilize this advantage.

\section{INCREMENT FORMULA}\label{sec:formula}

Given two controls $\bar u, u \in \mathcal U$, where $\bar u$ is an initial (reference) one, and $u \neq \bar u$ is the target one, we abbreviate by $\bar X$ and $X$ the flows of the vector fields $\bar V_t \doteq V_{\bar u(t)}$ and $V_t  \doteq V_{u(t)}$, respectively, and by $\mu: \, t \mapsto  \mu_t[u] = (X_t)\sharp \vartheta$ and $\bar \mu: \, t \mapsto  \mu_t[\bar u] = (\bar X_t)\sharp \vartheta$ the corresponding solutions to the Cauchy problem \eqref{conteq}, \eqref{mu-init}. 

Consider the increment 
\(
\Delta_{u} \mathcal I [\bar u] \doteq  \mathcal I[u] - \mathcal I[\bar u] \doteq \ell(\mu_T) - \ell(\bar \mu_T)
\)
of the cost functional. The base of our approach is the following result proved in Appendix \ref{appendix:formula-proof}.  
\begin{theorem}[Increment formula] Assume that $(A_{1})$--$(A_4)$ hold. Then, the following representation is valid:
\begin{align}
&\Delta_{u} \mathcal I [\bar u] = \label{formula}\\
&\int_0^T \!\d t \int D_\mu \ell^*\big|_{\left(\bar X_{t, T} \circ X_{t}\right)_\sharp\vartheta}\circ \bar X_{t, T} \, \bar J_{t}\,\left(V_t - \bar V_t\right)\d\mu_t.\notag
\end{align} 
Here, $\bar J$ denotes the solution of the linear problem \eqref{Jacobian} corresponding to $u=\bar u$ and $s=T$; $^*$ stands for the matrix transposition.\footnote{The term $D_\mu \ell^*\big|_{\left(\bar X_{t, T} \circ X_{t}\right)_\sharp\vartheta}\circ \bar X_{t, T}\,{\bar J_{t}}$ in \eqref{formula} is the gradient $\nabla_x \bar p_t^*$ of a characteristic solution $(t, x) \mapsto \bar p_t(x)$ to the dual transport equation of the form \eqref{dual-transp-loc} with $u=\bar u$ and the final condition $\bar p_T = \frac{\delta \ell}{\delta\mu}$.} %

\end{theorem}

Observe that formula \eqref{formula} represents the variation of $\mathcal I$ at the point $\bar u$ w.r.t. \emph{any other} admissible signal $u \in \mathcal U$; this formula is \emph{exact} (i.e. it does not contain any residuals). 

\begin{remark}
    The representation \eqref{formula} (and the consequent numeric method) can be literally adapted to the  case of distributed ensembles 
by replacing 
$(\bar V, \bar X, \bar J, X, V)$ with $(\bar V^\eta, \bar X^\eta, \bar J^\eta, X^\eta, V^\eta)$ and taking the expectation w.r.t. $\Xi$. 
\end{remark}

\subsection{Control Improvement 
}

The main consequence of the increment formula is the structure of controls of potential decrease from the reference point $\bar u$ provided by minimizers $w_t[\mu]$ in the problem
\begin{align}
&\min_{\upsilon \in U} \int \! D_\mu \ell^*\Big|_{\left(\bar X_{t, T} \circ X\right)_\sharp\vartheta} \!\!\!\!\!\circ \bar X_{t, T}\, \bar J_{t}\,  V_\upsilon  \, \d \mu\label{feedback}
\end{align}
viewed as $\mu$-feedback controls of the PDE \eqref{PDE}. Indeed, if $t\mapsto \check\mu_t$
is a well-defined solution to an initial value problem 
\eqref{PDE}, \eqref{mu-init}
with a backfed \emph{nonlocal} vector field $\check V_t \doteq  V_{w_t[\check\mu_t]}$,
and $u(t) \doteq w_t[\mu_t]$,
then, obviously,
\(
\Delta_{u} \mathcal I [\bar u] \leq 0.
\)
Thus, the cost of open-loop controls $u$ generated by the feedbacks \eqref{feedback} does not exceed (potentially, smaller than) the one of $\bar u$.

\subsection{Numeric Algorithm}\label{sec:alg}

A pitfall in the discussed control-update rule is due to the (generic) discontinuity of the map $x \mapsto \check V_t(x)$ that makes the Cauchy problem \eqref{flow} ill-posed. To resolve this issue, one can employ the classical semi-discrete Krasovskii-Subboting sampling scheme \cite{krasovskii2012control} with a time discretization (partition)
\(
\pi_{I}^N=\{0=t_0 < t_1 < \ldots < t_N = T\} \subset I.
\)

Let $u^k$, $k \in 0, 1, \ldots $, be given/computed. On the conceptual level, an iteration of the announced \underline{iterative method} consists of just three steps: 

i) integration of the 
ODE \eqref{flow} together with the linearized system 
\eqref{Jacobian}, for $\bar u=u^k$ and various initial conditions over 
some mesh 
\(
\pi_{\spt \vartheta}^M=\{y^k\}_{k=0}^M \subseteq \spt \vartheta\), to obtain $(X^k,J^k)$, 

ii) numeric solution of the PDE \eqref{PDE}, \eqref{mu-init}
backfed by \eqref{feedback} with $(\bar X, \bar J)=(X^k,J^k)$, to obtain $\mu^{k+1}$,
and 

iii) control update $u^{k+1} : = w_t[\mu^{k+1}_t]$.

Arguments similar to \cite[Appendix B]{SPP-2023} show that this iterative method 
converges in the residual of Pontryagin's maximum principle \cite{Bonnet2021AMT} for the convexified problem $(P)$ as \(\max\limits_{i}|t_{i} - t_{i-1}|+\max\limits_{k}\|y^{k}-y^{k-1}\| \to 0\) over $\pi^N_I \times \pi^M_{\spt \vartheta}$.\footnote{In \cite{Bonnet2021AMT}, the PMP is formulated for a $\mu$-linear problem with the functional of the form $\ell(\mu) = \int \varphi \d \mu$. Following the same line of reason, this result can be extended to a general $\C^1$ functional by an adequate modification of the transversality condition. 
}

\section{APPLICATION: BLOCH EQUATIONS}

We now apply the algorithm from \S~\ref{sec:alg} to a non-standard 
problem of designing composite pulses in a multi-population of nuclear spins, mentioned in the Introduction. Consider a family of  
Bloch equations, parameterized by  the (dimensionless)
resonance offset 
$\eta$. 
For simplicity, we focus on the non-dissipative case and rewrite the Bloch equations in spherical polar coordinates in the rotating frame \cite{Tah}:
 \begin{equation}
 \left(
 \begin{array}{c}\dot\theta \\ \dot\phi
 \end{array}
 \right) = V_u^\eta(\theta,\phi) \doteq u\left(
 \begin{array}{c}
 \cot \phi\cos \theta \\ \sin \theta\end{array}
 \right) - 
 \eta
 \left(\begin{array}{c}1\\0\end{array}
 \right)\label{Bloch}.
 \end{equation}
 Here, $\theta \in[0,2\pi]$ and $\phi\in [0,\pi]$ are the azimuthal and polar angles identifying the position on the Bloch sphere, $(x_1, x_2, x_3) \doteq (\cos\theta \sin\phi, \sin\theta \sin \phi, \cos\phi)$; 
control input $t \mapsto u(t)$ is the envelope of the actuating rf-field.\footnote{We restrict the control options to a single parameter representing the envelope of the exciting field, which essentially reduces the controllability and makes the resulting ensemble control problem much more challenging. 
}

\begin{remark}
It may be apt to stress that the Bloch equations
are not really of the quantum feature. These phenomenological ODEs describe the dynamics of an \emph{averaged} nuclear magnetization in a {macroscopic sample}, and are inapplicable to an individual  nuclear magnetic moment. In other words, 
each 
ODE \eqref{Bloch} already represents the dynamic ensemble. One can say that, in this example, we actually deal with an ``ensemble of ensembles''.  
\end{remark}

A canonical task in NMR experiments is to transfer the bulk magnetization vector from an equilibrium position (aligned with the static magnetic field) to the excited state $(\theta_T, \phi_T) = (0,\pi/2)$ (so-called $\pi/2$-transfer). In practice, the static field is inhomogeneous, which gives rise to probability distributions $\mu_0^\eta \in \mathcal P([0,2\pi]\times[0,\pi])$ in the initial values $(\theta_0, \phi_0)$, and leads to an optimal control problem of type \eqref{type-1}. We assume that $\mu_0^\eta$ are absolutely continuous with a common density function $\rho_0(\theta, \phi)$, 
and consider a more delicate performance criterion similar to \eqref{type2-problem} by incorporating a variance-like term and the energy cost \eqref{u-cost}. The resulting problem is adapted to the framework of distributed ensembles as follows:
\begin{equation}
\begin{array}{c}
\displaystyle\min \mathcal I[u] = \int_{-1}^1\ell(\mu_T^\eta)\d \Xi(\eta) + \frac{\alpha}{2}\int_0^T u^2(t) \d t,\\[0.3cm]
\ell(\mu) \doteq \,\displaystyle \int g(\cdot, \cdot, \theta_T, \phi_T) \d \mu + \frac{\beta}{2}\iint g \d(\mu\otimes\mu).\end{array}\label{cost}
\end{equation}
Here,
\(
g(\theta,\phi,\theta',\phi')=\frac{1}{2}[(\sin \theta - \sin \theta')^2 + (\sin \phi - \sin \phi')^2 + (\cos \theta - \cos \theta')^2 + (\cos \phi - \cos \phi')^2] \doteq 2-\cos(\theta-\theta') -\cos(\phi-\phi'),
\)
the integral $\int$ is computed over $[0,2\pi]\times[0,\pi]$, $\otimes$ denotes the tensor product of measures, and $\alpha, \beta>0$ are given parameters. 
To specify the feedback control \eqref{feedback}, we compute:
\(
D_\mu \ell(\mu)(\theta, \phi) =  D_{(\theta, \phi)}f(\theta, \phi) + \beta \displaystyle\int D_{(\theta, \phi)} g(\theta, \phi, \cdot) \d \mu.
\)

\begin{figure}[!t]
  \centerline{\includegraphics[width=0.42\textwidth]{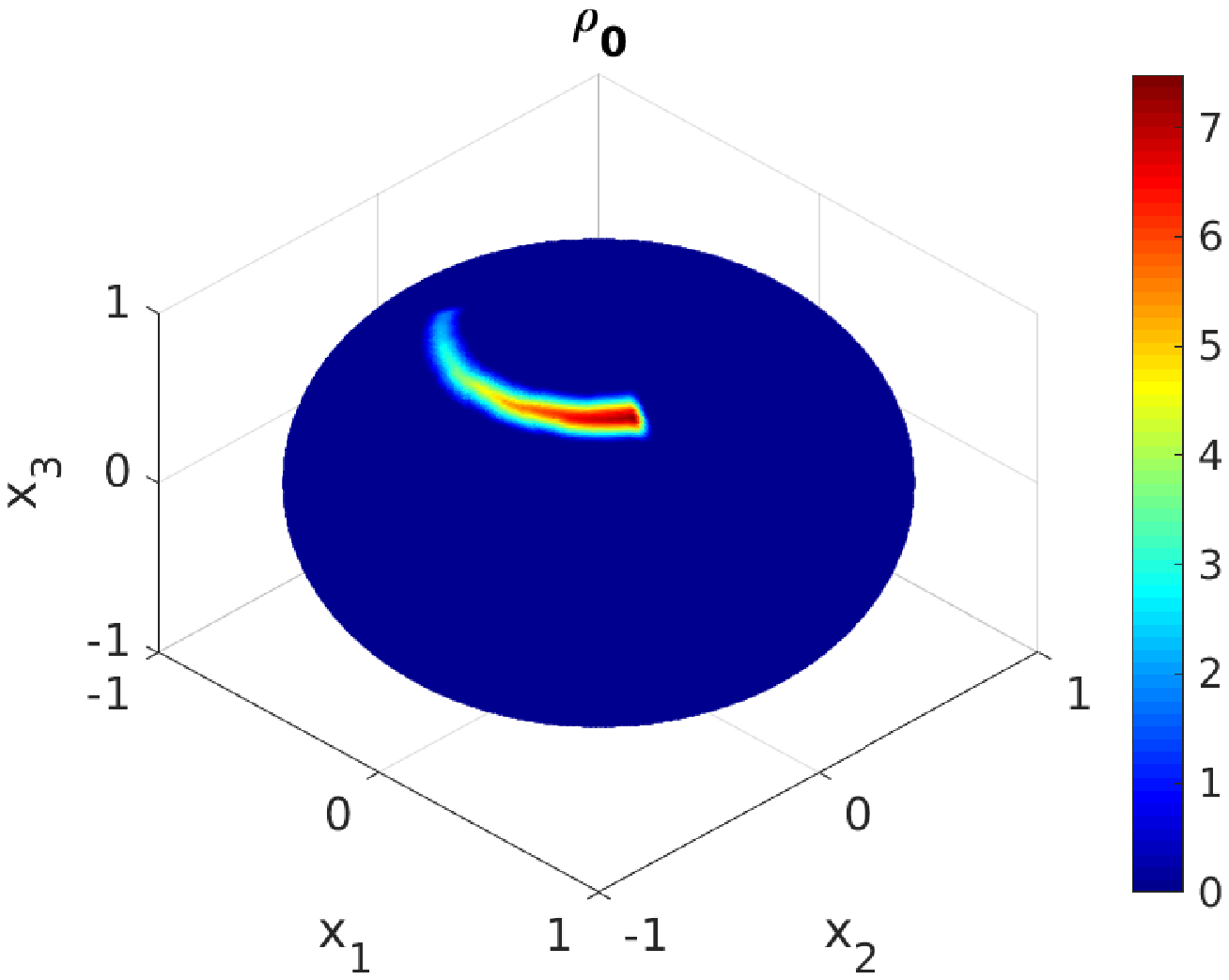}}  
  \centerline{\includegraphics[width=0.42\textwidth]{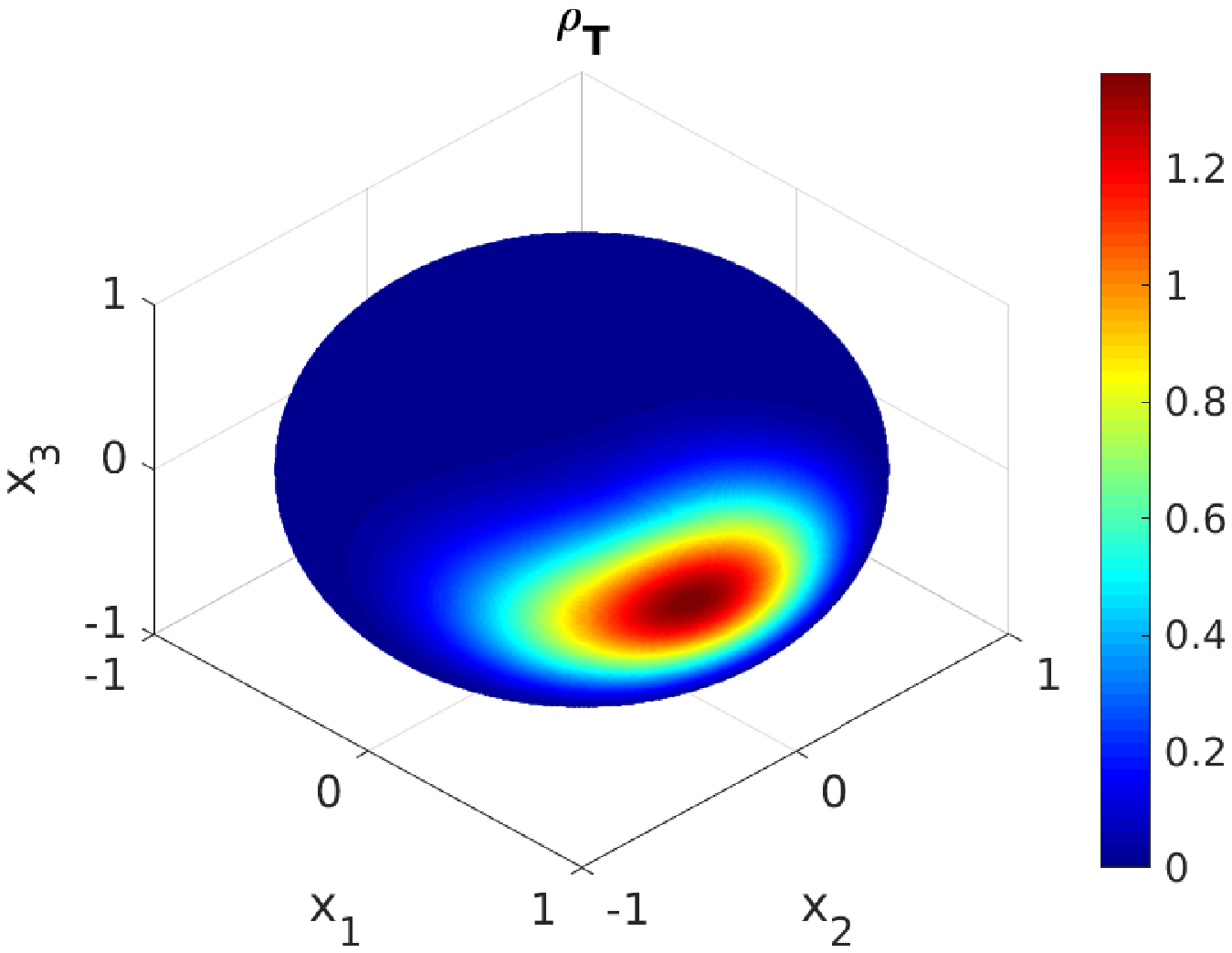}}  
  \centerline{\includegraphics[width=0.35\textwidth]{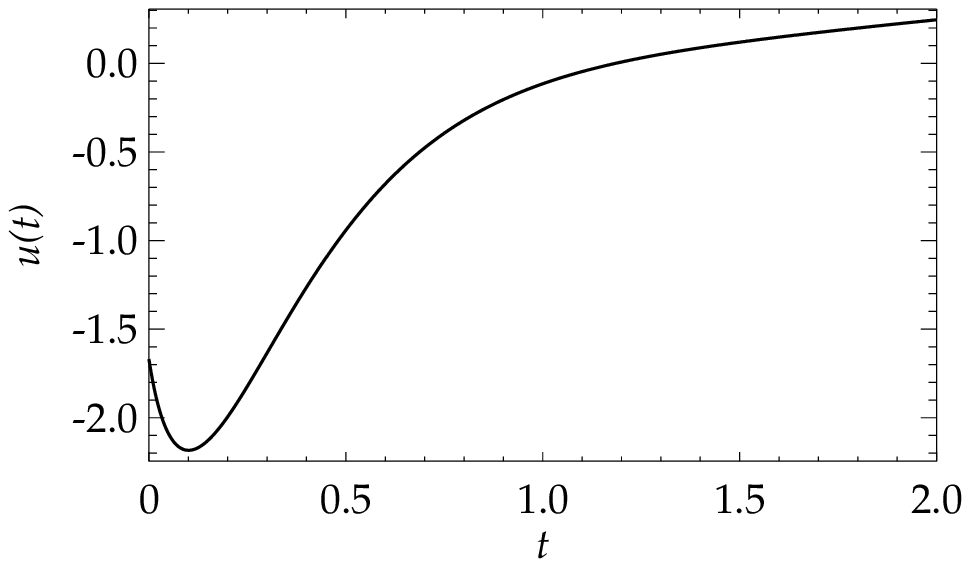}}  
  \caption{ Initial, for $t=0$, (top), final, for $t=T$, (top) density  $\rho_T$ and the corresponding control $u$ (bottom panel) after five iterations of the algorithm. 
  \label{fig:1}}
\end{figure}

We performed a numerical case study for the initial density $\rho_0$ (Fig.~\ref{fig:1}, top panel) on a uniform grid with a spacing of 0.01 for both angles,  
and for the distribution $\Xi$ chosen to be uniform on $[-0.55,-0.45]$ (due to the lack of space, the results are presented for the mean value $\eta=-0.5$). The standard Lax-Friedrichs numerical integration scheme was implemented (for integration in time from $0$ to $T=2$ with the constant time step $10^{-4}$). To exclude the singularity at the poles, the problem was solved for $\phi\in[0.05,0.095\pi]$, assuming the boundaries for $\theta$ to be periodic, and vanishing normal derivative for $\phi$.\footnote{For future simulations, in order to increase the computational efficiency of our codes and fix the pole problem, we plan to implement the pseudospectral methods using spherical harmonics. 
} The following values of the parameters in the cost function (\ref{cost}) were taken: $\alpha=0.25$ and $\beta=0.5.$ It turned out that the simulations are time consuming, thus, the code was parallelized for multiprocessor computers with shared memory. The simulations are also memory demanding~-- storage in memory of a large 
four-dimensional 
array (
a function of $(t,\theta,\phi,\eta)$) is required (about 150 GB for the parameter values described above).   

The initial control was taken to be constant, $u^0 \equiv 0.1$, with the  
cost 
$\mathcal I[u^0] \approx 0.88$. Computing the iterations of the proposed algorithms, the cost $\mathcal I$ 
was observed to decrease monotonically (as it is expected): $\approx$ 0.59, 0.49, 0.46, 0.44, 0.43 and then stagnating at a value $\approx 0.43$. Terminal density $\rho_T$ and the corresponding control $u$ computed after five iterations are shown in Fig.~\ref{fig:1} (middle and bottom panels, respectively).

\section{CONCLUSION}

We finally stress that the suggested nonlocal algorithm generates a $\mathcal I$-monotone control sequence, and typically takes a few (2-5) iterations to reach an acceptable solution. Free of 
any intrinsic parametric optimization, 
the method can be a ``lifeline'' for computationally demanding problems (like the presented one), where 
the direct approach, as well as different versions of the gradient descent, become prohibitively  expensive. 

Although the proposed approach has a fairly wide scope of application, there are 
significant restrictions. 
For instance, the concept of flat derivative (and, as a consequence, intrinsic derivative) does not apply to functionals $F: \, \mathcal P_c
\to \R$ that are undefined ($= +\infty$) for measures, singular w.r.t. some reference one (typically, $\mathcal L^n$); this makes it impossible to treat several useful performance criteria such as entropy functionals of the Kullback-Leibler type \cite{AmbrosioBook2007}. 

\bibliographystyle{IEEEtran}%
\bibliography{IEEEabrv,Staritsyn}

\def\cprime{$'$}
\begin{thebibliography}{10}
\providecommand{\url}[1]{#1}
\csname url@samestyle\endcsname
\providecommand{\newblock}{\relax}
\providecommand{\bibinfo}[2]{#2}
\providecommand{\BIBentrySTDinterwordspacing}{\spaceskip=0pt\relax}
\providecommand{\BIBentryALTinterwordstretchfactor}{4}
\providecommand{\BIBentryALTinterwordspacing}{\spaceskip=\fontdimen2\font plus
\BIBentryALTinterwordstretchfactor\fontdimen3\font minus
  \fontdimen4\font\relax}
\providecommand{\BIBforeignlanguage}[2]{{%
\expandafter\ifx\csname l@#1\endcsname\relax
\typeout{** WARNING: IEEEtran.bst: No hyphenation pattern has been}%
\typeout{** loaded for the language `#1'. Using the pattern for}%
\typeout{** the default language instead.}%
\else
\language=\csname l@#1\endcsname
\fi
#2}}
\providecommand{\BIBdecl}{\relax}
\BIBdecl

\bibitem{Li-Khaneja}
J.-S. Li and N.~Khaneja, ``Ensemble control of bloch equations,'' \emph{IEEE
  Transactions on Automatic Control}, vol.~54, no.~3, pp. 528--536, 2009.

\bibitem{Li-Khaneja2006}
------, ``Control of inhomogeneous quantum ensembles,'' \emph{Phys. Rev. A},
  vol.~73, p. 030302, Mar 2006.

\bibitem{SKINNER20038}
T.~E. Skinner, T.~O. Reiss, B.~Luy, N.~Khaneja, and S.~J. Glaser, ``Application
  of optimal control theory to the design of broadband excitation pulses for
  high-resolution {NMR},'' \emph{Journal of Magnetic Resonance}, vol. 163,
  no.~1, pp. 8--15, 2003.

\bibitem{Conolly}
S.~Conolly, D.~Nishimura, and A.~Macovski, ``Optimal control solutions to the
  magnetic resonance selective excitation problem,'' \emph{IEEE Transactions on
  Medical Imaging}, vol.~5, no.~2, pp. 106--115, 1986.

\bibitem{Bonnard-TAC}
B.~Bonnard, O.~Cots, S.~J. Glaser, M.~Lapert, D.~Sugny, and Y.~Zhang,
  ``Geometric optimal control of the contrast imaging problem in nuclear
  magnetic resonance,'' \emph{IEEE Transactions on Automatic Control}, vol.~57,
  no.~8, pp. 1957--1969, 2012.

\bibitem{Bonnard2016}
B.~Bonnard, A.~Jacquemard, and J.~Rouot, ``Optimal control of an ensemble of
  {Bloch} equations with applications in mri,'' in \emph{2016 IEEE 55th
  Conference on Decision and Control (CDC)}, 2016, pp. 1608--1613.

\bibitem{WANG2018306}
S.~Wang and J.-S. Li, ``Free-endpoint optimal control of inhomogeneous bilinear
  ensemble systems,'' \emph{Automatica}, vol.~95, pp. 306--315, 2018.

\bibitem{CHEN2020109057}
X.~Chen, ``Ensemble observability of bloch equations with unknown population
  density,'' \emph{Automatica}, vol. 119, p. 109057, 2020.

\bibitem{AmbrosioBook2007}
L.~Ambrosio and G.~Savar\'e, ``Gradient flows of probability measures,'' in
  \emph{Handbook of differential equations: evolutionary equations. {V}ol.
  {III}}, ser. Handb. Differ. Equ.\hskip 1em plus 0.5em minus 0.4em\relax
  Elsevier/North-Holland, Amsterdam, 2007, pp. 1--136.

\bibitem{Cavagnari2018}
G.~Cavagnari, A.~Marigonda, K.~T. Nguyen, and F.~S. Priuli, ``Generalized
  control systems in the space of probability measures,'' \emph{Set-Valued and
  Variational Analysis}, vol.~26, no.~3, pp. 663--691, 2018.

\bibitem{Bonnet2021AMT}
\BIBentryALTinterwordspacing
B.~Bonnet, C.~Cipriani, M.~Fornasier, and H.-L. Huang, ``A measure theoretical
  approach to the mean-field maximum principle for training neurodes,'' 2021.
  [Online]. Available: \url{https://arxiv.org/abs/2107.08707}
\BIBentrySTDinterwordspacing

\bibitem{SChPP-2022}
M.~Staritsyn, N.~Pogodaev, R.~Chertovskih, and F.~L. Pereira, ``Feedback
  maximum principle for ensemble control of local continuity equations: An
  application to supervised machine learning,'' \emph{IEEE Control Systems
  Letters}, vol.~6, pp. 1046--1051, 2022.

\bibitem{SPP-2023}
M.~Staritsyn, N.~Pogodaev, and F.~L. Pereira, ``Linear-quadratic problems of
  optimal control in the space of probabilities,'' \emph{IEEE Control Systems
  Letters}, vol.~6, pp. 3271--3276, 2022.

\bibitem{CardMaster2019}
P.~Cardaliaguet, F.~Delarue, J.-M. Lasry, and P.-L. Lions, \emph{The Master
  Equation and the Convergence Problem in Mean Field Games}, ser. Ann.
  {{Math}}. {{Stud}}.\hskip 1em plus 0.5em minus 0.4em\relax {Princeton, NJ:
  Princeton University Press}, 2019, vol. 201.

\bibitem{POGODAEV20203585}
N.~Pogodaev and M.~Staritsyn, ``Impulsive control of nonlocal transport
  equations,'' \emph{Journal of Differential Equations}, vol. 269, no.~4, pp.
  3585--3623, 2020.

\bibitem{Li-TACON-2013}
J.-S. Li, I.~Dasanayake, and J.~Ruths, ``Control and synchronization of neuron
  ensembles,'' \emph{IEEE Transactions on Automatic Control}, vol.~58, no.~8,
  pp. 1919--1930, 2013.

\bibitem{ZUAZUA20143077}
E.~Zuazua, ``Averaged control,'' \emph{Automatica}, vol.~50, no.~12, pp.
  3077--3087, 2014.

\bibitem{Bressan-PiccoliBOOK}
A.~Bressan and B.~Piccoli, \emph{Introduction to the mathematical theory of
  control}, ser. AIMS Series on Applied Mathematics.\hskip 1em plus 0.5em minus
  0.4em\relax American Institute of Mathematical Sciences (AIMS), Springfield,
  MO, 2007, vol.~2.

\bibitem{Agrachev-Book}
A.~Agrachev and Y.~Sachkov, \emph{Control Theory from the Geometric Viewpoint},
  ser. Control theory and optimization.\hskip 1em plus 0.5em minus 0.4em\relax
  Springer, 2004.

\bibitem{Bonnet2021}
B.~Bonnet and H.~Frankowska, ``Necessary optimality conditions for optimal
  control problems in {Wasserstein} spaces,'' \emph{Applied Mathematics {\&}
  Optimization}, vol.~84, no.~2, pp. 1281--1330, Dec 2021.

\bibitem{krasovskii2012control}
A.~Krasovskii and N.~Krasovskii, \emph{Control Under Lack of Information}, ser.
  Systems \& Control: Foundations \& Applications.\hskip 1em plus 0.5em minus
  0.4em\relax Birkh{\"a}user Boston, 2012.

\bibitem{Tah}
B.~Tahayori, L.~Johnston, I.~Mareels, and P.~Farrell, ``Novel insight into
  magnetic resonance through a spherical coordinate framework for the bloch
  equation,'' in \emph{Progress in Biomedical Optics and Imaging - Proceedings
  of SPIE}, vol. 7258, 02 2009.

\end{thebibliography}

\appendix

\subsection{Proof of Theorem 1}\label{appendix:formula-proof}

Let $t\mapsto \mu_{t}$ and $t\mapsto \bar \mu_{t}$ denote the weak solutions of the PDE 
\eqref{conteq} with initial condition $\mu_{0} = \vartheta$, corresponding to control inputs $u$ and $\bar u$, respectively. Recall that
\(
\mu_t = (X_{0,t})_\sharp \vartheta\), \(\bar \mu_t = (\bar X_{0,t})_\sharp \vartheta, 
\)
and
\(
X_{s,s} = \bar X_{s,s} = \id\) \(\forall s \in \R,
\)
where $X$ and $\bar X$ are the corresponding characteristic flows.

Denote $\mathcal F_t \doteq \bar X_{t, T} \circ X_{0,t}$. Since the map $x \mapsto \mathcal F_t(x)$ is a composition of two bijections, it is invertible. Standard arguments from the ODE theory imply that under assumptions $(A_1)$--$(A_3)$, the maps $t \mapsto \bar X_{t, T}(x)$ and $t \mapsto X_{0, t}(x)$ are Lipschitz on $I$, for any $x \in \R^n$. Hence, for any $x \in \R^n$, the function $t \mapsto \mathcal F_t(x)$ is absolutely continuous on $I$ as a composition of Lipschitz maps; in particular it is  $\mathcal L^1$-a.e. differentiable: 
\(
\mathcal F_{t+\lambda} = \mathcal F_{t} + \lambda \, \partial_t \mathcal F_{t} + o(\lambda) \approx (\id + \lambda \, \partial_t \mathcal F_{t} \circ \mathcal F_{t}^{-1}) \circ \mathcal F_{t}.
\)
Assumption $(A_4)$ guarantees that, for any Borel measurable, locally bounded map \( \varphi\colon \R^{n}\to \R^{n} \), the function \( \lambda\mapsto \ell\left((\id+\lambda\mathcal F)_{\sharp}\mu\right) \) is differentiable at zero, and
\(\displaystyle\frac{d}{d\lambda}\Big|_{\lambda=0}\ell\left((\id+\lambda f)_{\sharp}\mu\right) = \int D_{\mu}\ell(\mu)\cdot f\d\mu\),
where $D_{\mu}\ell$ stands for the intrinsic derivative. Thus,
\begin{align}
\partial_t \ & \ \ell\left(\left(\mathcal F_t\right)_\sharp \vartheta\right) =  \notag \\ 
= \ & \ \frac{d}{d \lambda}\Big|_{\lambda=0} \ell\left(\left(\mathcal F_{t+\lambda}\right)_\sharp \vartheta\right) \notag \\
= \ & \ \frac{d}{d \lambda}\Big|_{\lambda=0} \ell\left(\left(\id + \lambda \, \partial_t \mathcal F_{t}\circ \mathcal F_{t}^{-1}\right)_\sharp \big(\left(\mathcal F_t\right)_\sharp \vartheta\big)\right)\notag\\
= \  & \ \int D_\mu \ell\left(\left(\mathcal F_t\right)_\sharp \vartheta\right) \cdot \partial_t \mathcal F_t \circ \mathcal F_{t}^{-1} \d{\big(\left(\mathcal F_t\right)_\sharp \vartheta\big)}\notag\\=
\ & \ \int D_\mu 
\ell\big|_{\left(\bar X_{t, T} \circ X_{t}\right)_\sharp\vartheta}
\circ \left(\bar X_{t, T} \circ X_{0,t}\right) 
\cdot \notag\\
& \qquad \partial_t (\bar X_{t, T} \circ X_{0,t}) \d \vartheta
.\label{der-nonlocal-formula}
\end{align}

In the last expression, the partial derivative in $t$ is represented by the chain rule as
\[
\partial_t (\bar X_{t, T} \circ X_{0,t}) =  \left[\partial_\tau \bar X_{t, T} \circ X_{0,\tau} + \partial_\tau \bar X_{\tau, T} \circ X_{0,t}\right]\big|_{\tau =t},
\]
where
\(
\partial_\tau\big|_{\tau =t} \bar X_{t, T} \circ X_{0,\tau} = \left(D_x \bar X_{t, T}  \, V_t\right)\circ X_{0,t}
\doteq \left(\bar J_{t, T}  \, V_t\right)\circ X_{0,t}\)
by direct computation, and
\(
\partial_\tau\big|_{\tau =t} \bar X_{\tau, T} = -\bar J_{t,T}  \, \bar V_t
\)
by \eqref{inv-flow-derivative}. 
Plugging these expressions to \eqref{der-nonlocal-formula}, we obtain
\begin{align}
\partial_t 
\ell\left(\left(\mathcal F_t\right)_\sharp \vartheta\right) 
\doteq & \ \int \Big[D_\mu \ell^*\big|_{\left(\bar X_{t, T} \circ X_{t}\right)_\sharp\vartheta}
\circ \bar X_{t, T} 
\notag\\& \qquad 
\bar J_{t, T} (V_t-\bar V_t)\Big]  \circ X_{0,t} \d \vartheta.\label{ell-proof}
\end{align}

Now, the cost increment 
is represented as follows:
\begin{align*}
    \Delta_{u}\mathcal I[\bar u] & \doteq \ell(\mu_T) - \
    \ell(\bar \mu_T)\\
    & = \ell\left((\bar X_{T,T} \circ X_{0,T})_\sharp \vartheta\right) - \ell\left((\bar X_{T,T} \circ \bar X_{0,T})_\sharp \vartheta\right)\\
    & - \underbrace{\left[\ell\left((\bar X_{0,T} \circ X_{0,0})_\sharp \vartheta\right) - \ell\left((\bar X_{0,T} \circ \bar X_{0,0})_\sharp \vartheta\right)\right]}_\text{$ \equiv 0$}\\[-0.2cm]
    = \int_0^T & \!\! \partial_t \big[\ell\left((\bar X_{t,T} \circ X_{0,t})_\sharp \vartheta\right) - \ell\left((\bar X_{t,T} \circ \bar X_{0,t})_\sharp \vartheta\right)\big]\d t.
\end{align*}
By the semigroup property, \(\bar X_{t, T}\circ\bar X_{0,t} = \bar X_{0,T}\), which implies that the second term under the sign of the time derivative in the latter expression is, in fact, independent of $t$, and therefore, $\Delta_{u}\mathcal I[\bar u]$ equals
\begin{align*}
\int_0^T \partial_t \ell\left(\left(\bar X_{t, T} \circ X_{0,t}\right)_\sharp \vartheta\right)\d t =\int_0^T \partial_t\ell\left(\left(\mathcal F_t\right)_\sharp \vartheta\right)\d t.\label{incr-proof}
\end{align*}
To complete the proof, it remains to combine the latter expression with \eqref{ell-proof} and use the representation formula \(
\mu_t = (X_{0,t})_\sharp \vartheta
\).

\end{document}